\documentclass{amsart}

\usepackage{amssymb}

\allowdisplaybreaks

\newtheorem{theorem}{Theorem}[section]
\newtheorem{lemma}[theorem]{Lemma}

\newtheorem{corollary}[theorem]{Corollary}
\newtheorem{claim}[theorem]{Claim}
\newtheorem{remark}[theorem]{Remark}

\theoremstyle{definition}

\newcommand{\dem}{\noindent {\bf Proof. }}
\newcommand{\sketch}{\noindent {\bf Sketch of the proof. }}
\newcommand{\fin}{\hspace*{\fill} $\square$ \vskip0.2cm}

\newcommand{\rad}{\mathrm{r}}

\newcommand{\Rad}{\mathcal{R}}
\newcommand{\Word}{\mathcal{W}}
\newcommand{\Max}{\mathcal{K}}

\newcommand{\Free}{\mathbf{F}}

\def\astop{\mathop{*}}

\begin{document}

\title[Khintchine type inequalities associated with Free groups]
{Non-commutative Khintchine type inequalities \\ associated with
free groups}

\author[Parcet and Pisier]
{Javier Parcet$^*$ and Gilles Pisier$^{\dag}$}

\date{}

\begin{abstract}
Let $\Free_n$ denote the free group with $n$ generators $g_1, g_2,
\ldots, g_n$. Let $\lambda$ stand for the left regular
representation of $\Free_n$ and let $\tau$ be the standard trace
associated to $\lambda$. Given any positive integer $d$, we study
the operator space structure of the subspace $\Word_p(n,d)$ of
$L_p(\tau)$ generated by the family of operators
$\lambda(g_{i_1}g_{i_2} \cdots g_{i_d})$ with $1 \le i_k \le n$.
Moreover, our description of this operator space holds up to a
constant which does not depend on $n$ or $p$, so that our result
remains valid for infinitely many generators. We also consider the
subspace of $L_p(\tau)$ generated by the image under $\lambda$ of
the set of reduced words of length $d$. Our result extends to any
exponent $1 \le p \le \infty$ a previous result of Buchholz for
the space $\Word_{\infty}(n,d)$. The main application is a certain
interpolation theorem, valid for any degree $d$ (extending a
result of the second author restricted to $d=1$). In the simplest
case $d=2$, our theorem can be stated as follows: consider the
space $\mathcal{K}_p$ formed of all block matrices $a=(a_{ij})$
with entries in the Schatten class $S_p$, such that $a$ is in
$S_p$ relative to $\ell_2 \otimes \ell_2$ and moreover such that
$(\sum_{ij} a_{ij}^* a_{ij} )^{1/2}$ and $(\sum_{ij} a_{ij}
a_{ij}^*)^{1/2}$ both belong to $S_p$. We equip $\mathcal{K}_p$
with the maximum of the three corresponding norms. Then, for $2
\le p \le \infty$ we have $\mathcal{K}_p \simeq (\mathcal{K}_2,
\mathcal{K}_\infty)_\theta$ with $1/p = (1-\theta)/2$.
\end{abstract}

\maketitle

\section*{Introduction}

Let $\Rad_p(n)$ be the subspace of $L_p[0,1]$ generated by the
classical Rademacher functions $\rad_1, \rad_2, \ldots, \rad_n$.
As is well-known, for any exponent $1 \le p < \infty$, the
classical Khintchine inequalities provide a linear isomorphism
between $\Rad_p(n)$ and $\ell_2(n)$ with constants independent of
$n$. However, to describe the operator space structure of
$\Rad_p(n)$ we need the so-called non-commutative Khintchine
inequalities, introduced by F. Lust-Piquard in \cite{L} and
extended in \cite{LP} to the case $p=1$, see also \cite{Bu2} for
an analysis of the optimal constants.

\vskip3pt

To describe the non-commutative Khintchine inequalities, let us
consider a family $x_1, x_2, \ldots, x_n$ of elements in the
Schatten class $S_p$. Then we have the following equivalences of
norms for $1 \le p \le 2$ $$\Big\| \sum_{k=1}^n x_k \rad_k
\Big\|_{L_p([0,1];S_p)} \simeq \, \inf_{x_k = y_k + z_k} \left\{
\Big\| \Big( \sum_{k=1}^n y_k y_k^* \Big)^{1/2} \Big\|_{S_p} +
\Big\| \Big( \sum_{k=1}^n z_k^* z_k \Big)^{1/2} \Big\|_{S_p}
\right\},$$ while for $2 \le p < \infty$ we have $$\Big\|
\sum_{k=1}^n x_k \rad_k \Big\|_{L_p([0,1];S_p)} \simeq \, \max
\left\{ \Big\| \Big( \sum_{k=1}^n x_k x_k^* \Big)^{1/2}
\Big\|_{S_p} , \Big\| \Big( \sum_{k=1}^n x_k^* x_k \Big)^{1/2}
\Big\|_{S_p} \right\}.$$ Again the constants do not depend on $n$.
According to the purposes of this paper, it will be more
convenient to rewrite these inequalities in terms of the spaces
$R_p^n$ and $C_p^n$. Let us consider the Schatten class $S_p^n$
over the $n \times n$ matrices and let us denote its natural basis
by $$\Big\{ e_{ij} \, | \ 1 \le i,j \le n \Big\}.$$ We define
$R_p^n$ to be the subspace of $S_p^n$ generated by $e_{11},
e_{12}, \ldots, e_{1n}$ while $C_p^n$ will be the subspace
generated by $e_{11}, e_{21}, \ldots, e_{n1}$. That is, $R_p^n$
and $C_p^n$ can be regarded as the row and column subspaces of
$S_p^n$ respectively. Note that both subspaces have a natural
operator space structure inherited from $S_p^n$. In terms of the
spaces $R_p^n$ and $C_p^n$, the non-commutative Khintchine
inequalities can be rephrased by saying that $\Rad_p(n)$ is
completely isomorphic to $R_p^n + C_p^n$ whenever $1 \le p \le 2$
and $\Rad_p(n)$ is completely isomorphic to $R_p^n \cap C_p^n$ for
$2 \le p < \infty$. More concretely, when $1 \le p \le 2$ it
follows that
\begin{eqnarray*}
\lefteqn{\Big\| \sum_{k=1}^n x_k \rad_k \Big\|_{L_p([0,1];S_p)}}
\\ & \simeq & \inf_{x_k = y_k + z_k} \left\{ \Big\| \sum_{k=1}^n y_k
\otimes e_{1k} \Big\|_{S_p(\ell_2 \otimes \ell_2)} + \Big\|
\sum_{k=1}^n z_k \otimes e_{k1} \Big\|_{S_p(\ell_2 \otimes
\ell_2)} \right\},
\end{eqnarray*}
while for $2 \le p \le \infty$ we have
\begin{eqnarray*}
\lefteqn{\Big\| \sum_{k=1}^n x_k \rad_k \Big\|_{L_p([0,1];S_p)}}
\\ & \simeq & \max \left\{ \Big\| \sum_{k=1}^n x_k \otimes e_{1k}
\Big\|_{S_p(\ell_2 \otimes \ell_2)}, \Big\| \sum_{k=1}^n x_k
\otimes e_{k1} \Big\|_{S_p(\ell_2 \otimes \ell_2)} \right\}.
\end{eqnarray*}

Now let $\Free_n$ be the free group with $n$ generators $g_1,g_2,
\ldots, g_n$. If $\lambda$ denotes the left regular representation
of $\Free_n$, the family of operators $\lambda(g_1), \lambda(g_2),
\ldots, \lambda(g_n)$ appear as the free analog of the sequence
$\rad_1, \rad_2, \ldots, \rad_n$ in this framework. Namely, let us
consider the standard trace $\tau$ on $C_{\lambda}^*(\Free_n)$.
Then, by \cite{HP} it turns out that the subspace $\Word_p(n)$ of
$L_p(\tau)$ generated by the operators $\lambda(g_1),
\lambda(g_2), \ldots, \lambda(g_n)$ is completely isomorphic to
$\Rad_p(n)$ (with constants independent of $n$) for $1 \le p <
\infty$. Moreover, for $p = \infty$ we have
\begin{eqnarray*}
\lefteqn{\Big\| \sum_{k=1}^n a_k \otimes \lambda(g_k)
\Big\|_{L_{\infty}(\mathrm{tr} \otimes \tau)}} \\ & \simeq & \max
\left\{ \Big\| \sum_{k=1}^n a_k \otimes e_{1k}
\Big\|_{S_{\infty}(\ell_2 \otimes \ell_2)} , \Big\| \sum_{k=1}^n
a_k \otimes e_{k1} \Big\|_{S_{\infty}(\ell_2 \otimes \ell_2)}
\right\},
\end{eqnarray*}
where $R_n$ and $C_n$ are the usual expressions for $R_{\infty}^n$
and $C_{\infty}^n$. In other words, we also have
$\mathcal{W}_{\infty}(n) \simeq R_n \cap C_n$ completely
isomorphically. In this paper we shall generalize the mentioned
complete isomorphism for $\Word_p(n)$, where the results above
appear as the case of degree one. More concretely, given any
positive integer $d$, we shall consider (operator-valued)
homogeneous polynomials of degree $d$ in the variables
$$\lambda(g_1), \lambda(g_2), \ldots, \lambda(g_n).$$

\vskip3pt

Let $\Word_p(n,d)$ be the subspace of $L_p(\tau)$ generated by the
operators $\lambda(g_{i_1}g_{i_2} \cdots g_{i_d})$, where $1 \le
i_k \le n$ for $1 \le k \le d$. The aim of this paper is to
describe the operator space structure of $\Word_p(n,d)$ for any
value of $d$. To that aim, we consider an auxiliary
non-commutative $L_p$ space $L_p(\varphi)$ equipped with a
faithful normal semi-finite trace $\varphi$. Then, describing the
operator space structure of $\Word_p(n,d)$ becomes equivalent to
describing the norm of $$\sum_{i_1, i_2, \ldots, i_d = 1}^n
a_{i_1i_2 \cdots i_d} \otimes \lambda(g_{i_1}g_{i_2} \cdots
g_{i_d})$$ in $L_p(\varphi \otimes \tau)$ up to constants not
depending on $n$ or $p$. Let $\mathcal{A}$ be the family of
operators $\big\{a_{i_1i_2 \cdots i_d} \, | \ 1 \le i_k \le n
\big\}$. This family can be regarded as an element of a
matrix-valued $L_p(\varphi)$-space in several ways. Namely, given
$0 \le k \le d$, we can construct a matrix $\mathcal{A}_k$ with
entries in $\mathcal{A}$ by taking the first $k$ indices $i_1,
i_2, \ldots, i_k$ as the row index and the last $d-k$ indices
$i_{k+1}, i_{k+2}, \ldots, i_d$ as the column index. In other
words, we consider the matrix $$\mathcal{A}_k = \Big( a_{(i_1
\cdots i_k),(i_{k+1} \cdots i_d)} \Big).$$ In particular, we see
$\mathcal{A}_k$ as an element of the Haagerup tensor product
$C_p^{n^k} \otimes_h R_p^{n^{d-k}}$ with values in $L_p(\varphi)$.
This allows us to define the following family of spaces
$$\begin{array}{lcl} \displaystyle \Max_p(n,d) = \sum_{k=0}^d
C_p^{n^k} \otimes_h R_p^{n^{d-k}} & \qquad & \mbox{for} \ 1 \le p
\le 2, \\ \displaystyle \Max_p(n,d) = \bigcap_{k=0}^d C_p^{n^k}
\otimes_h R_p^{n^{d-k}} & \qquad & \mbox{for} \ 2 \le p \le
\infty.
\end{array}$$

In this paper we shall prove that $\Word_p(n,d)$ and $\Max_p(n,d)$
are completely isomorphic operator spaces with constants depending
only on the degree $d$. More concretely, the following result
holds.

\vskip3pt

\noindent \textbf{Theorem.} \emph{$\Max_p(n,d)$ and $\Word_p(n,d)$
are completely isomorphic for $1 \le p \le \infty$ and there
exists an absolute constant $\mathrm{c}_d$ depending only on $d$
such that the following inequalities hold $$\frac{1}{\mathrm{c}_d}
\, \|\mathcal{A}\|_{L_p(\varphi; \Max_p(n,d))} \le \Big\|
\sum_{i_1 \ldots i_d=1}^n a_{i_1 \cdots i_d} \otimes
\lambda(g_{i_1} \cdots g_{i_d}) \Big\|_{L_p(\varphi \otimes \tau)}
\le \mathrm{c}_d \, \|\mathcal{A}\|_{L_p(\varphi; \Max_p(n,d))}.$$
Moreover, the natural projection $\mathrm{P}: L_p(\tau)
\rightarrow \Word_p(n,d)$ is c.b. with $\|\mathrm{P}\|_{cb} \le
\mathrm{c}_d$.}

\vskip3pt

The case $p=\infty$ of this result is the content of Buchholz's
paper \cite{Bu1}. However, a slightly different definition of
$\Word_{\infty}(n,d)$ was used there. Indeed, Buchholz takes the
subspace of $L_{\infty}(\tau)$ generated by the image under
$\lambda$ of the set of reduced words of length $d$. That is,
inverses of generators are also allowed to appear. As we shall
see, our arguments apply in this context for any $1 \le p \le
\infty$ and we will also provide the proof. In Section
\ref{Section-Iteration} we set some preliminary results while
Section \ref{Section-Main} is devoted to the proof of our main
result and its analog for reduced words of length $d$.

\section{Preliminary results}
\label{Section-Iteration}

One of the key points in the proof of our main result lies in the
iteration of the non-commutative Khintchine inequality applied to
the operator space $\Word_p(n)^{\otimes d}$. This process was
already pointed in \cite[Section 9.8]{P2}, but we provide here an
explicit description of the resulting inequalities. After that, we
state an $L_p$-valued version of Fell's absorption principle that
will be also needed in the proof.

\subsection{Iterations of the Khintchine inequality}

Let $\Free_n$ be the free group with $n$ generators $g_1,g_2,
\ldots, g_n$. If $(\delta_t)_{t \in \Free_n}$ denotes the natural
basis of $\ell_2(\Free_n)$, the left regular representation
$\lambda$ of $\Free_n$ is defined by the relation $\lambda(t_1)
\delta_{t_2} = \delta_{t_1t_2}$. The reduced $C^*$-algebra
$C_{\lambda}^*(\Free_n)$ is defined as the $C^*$-algebra generated
in $\mathcal{B}(\ell_2(\Free_n))$ by the operators $\lambda(t)$
when $t$ runs over $\Free_n$. Let us denote by $\tau$ the standard
trace on $C_{\lambda}^*(\Free_n)$ defined by $\tau(x) = \langle x
\delta_e, \delta_e \rangle$, where $e$ denotes the identity
element of $\Free_n$. Then, we construct the non-commutative $L_p$
space $L_p(\tau)$ in the usual way and consider the subspace
$\Word_p(n)$ of $L_p(\tau)$ generated by the operators
$\lambda(g_1), \lambda(g_2), \ldots, \lambda(g_n)$. The aim of
this section is to describe the operator space structure of
$\Word_p(n)^{\otimes d}$ as a subspace of $L_p(\tau^{\otimes d})$
for the exponents $2 \le p \le \infty$. As it was pointed out in
\cite{P2}, the case $1 \le p \le 2$ follows easily by duality.
However, we shall not write the explicit inequalities in that case
since we are not using them and the notation is considerably more
complicated.

\vskip3pt

The following result can be regarded as a particular case of our
main result for homogeneous polynomials of degree $1$. It was
proved in \cite{HP} for $p = \infty$ while the proof for $2 \le p
< \infty$ can be found in Corollary 9.7.2 of \cite{P2}. We notice
that its proof uses the fact that $R_p^n \cap C_p^n$ is an
interpolation family for $2 \le p \le \infty$.

\begin{lemma} \label{Lemma-Iteration}
The following equivalence of norms holds for $2 \le p \le \infty$,
$$\Big\| \sum_{k=1}^n a_k \otimes \lambda(g_k) \Big\|_{L_p(\varphi
\otimes \tau)} \simeq \max \left\{ \Big\| \sum_{k=1}^n a_k \otimes
e_{1k} \Big\|_{L_p(\varphi;R_p^n)} , \Big\| \sum_{k=1}^n a_k
\otimes e_{k1} \Big\|_{L_p(\varphi;C_p^n)} \right\} \! .$$ In
fact, the linear map $u: R_p^n \cap C_p^n \rightarrow \Word_p(n)$
defined by $$u(e_{1k} \oplus e_{k1}) = \lambda(g_k),$$ is a
complete isomorphism with $\|u\|_{cb} \le 2$ and completely
contractive inverse. On the other hand, the canonical projection
$\mathrm{P}: L_p(\tau) \rightarrow \Word_p(n)$ satisfies
$\|\mathrm{P}\|_{cb} \le 2$.
\end{lemma}

Let us consider the group product $\mathrm{G}_d = \Free_n \times
\Free_n \times \cdots \times \Free_n$ with $d$ factors. The left
regular representation $\lambda_d$ of $\mathrm{G}_d$ has the form
$$\lambda_d(t_1,t_2, \ldots, t_d) = \lambda(t_1) \otimes
\lambda(t_2) \otimes \cdots \otimes \lambda(t_d),$$ where
$\lambda$ still denotes the left regular representation of
$\Free_n$. In particular, the reduced $C^*$-algebra
$C_{\lambda_d}^*(\mathrm{G}_d)$ is endowed with the trace $\tau_d
= \tau \otimes \tau \otimes \cdots \otimes \tau$ with $d$ factors.
This allows us to consider the non-commutative space $L_p(\tau_d)$
for any $1 \le p \le \infty$. Then we define the space
$\Word_p(n)^{\otimes d}$ to be the subspace of $L_p(\tau_d)$
generated by the family of operators $$\lambda(g_{i_1}) \otimes
\lambda(g_{i_2}) \otimes \cdots \otimes \lambda(g_{i_d}).$$ If we
apply repeatedly Lemma \ref{Lemma-Iteration} to the sum
$$\mathcal{S}_d(a) = \sum_{i_1, \ldots, i_d=1}^n a_{i_1 i_2 \cdots
i_d} \otimes \lambda(g_{i_1}) \otimes \lambda(g_{i_2}) \otimes
\cdots \otimes \lambda(g_{i_d}) \in L_p(\varphi \otimes \tau_d),$$
then we easily get $$\|\mathcal{S}_d(a)\|_{L_p(\varphi \otimes
\tau_d)} \le 2^d \max \left\{ \Big\| \sum_{i_1, \ldots, i_d =1}^n
a_{i_1 \cdots i_d} \otimes \xi_1(i_1) \otimes \cdots \otimes
\xi_d(i_d) \Big\|_{S_p^{n^d}(L_p(\varphi))} \ \right\},$$ where
the maximum runs over all possible ways to choose the functions
$\xi_1, \xi_2, \ldots, \xi_d$ among $\xi_k(\cdot) = e_{\cdot 1}$
and $\xi_k(\cdot) = e_{1 \cdot}$. That is, each function $\xi_k$
can take values either in the space $R_p^n$ or in the space
$C_p^n$. For a given selection of $\xi_1, \xi_2, \ldots, \xi_d$ we
split up these functions into two sets, one made up of the
functions taking values in $R_p^n$ and the other taking values in
$C_p^n$. More concretely, let us consider the sets
$$\begin{array}{l} \mathrm{R}_{\xi} = \big\{ k \, | \ \xi_k(i) =
e_{1i} \big\}, \\ \mathrm{C}_{\xi} = \big\{ k \, | \ \xi_k(i) =
e_{i1} \big\}. \end{array}$$ Then, if $\mathrm{C}_{\xi}$ has $s$
elements, the sum $$\sum_{i_1, \ldots, i_d =1}^n a_{i_1 \cdots
i_d} \otimes \xi_1(i_1) \otimes \cdots \otimes \xi_d(i_d)$$ can be
regarded as a $n^s \times n^{d-s}$ matrix with entries in
$L_p(\varphi)$. Now we introduce a simpler notation already
employed in \cite{HP}. Let $[m]$ be an abbreviation for the set
$\{1,2, \ldots, m\}$. Then, if $\mathbb{P}_d(2)$ denotes the set
of partitions $(\alpha, \beta)$ of $[d]$ into two disjoint subsets
$\alpha$ and $\beta$, we denote by $$\pi_{\alpha}: [n]^d
\rightarrow [n]^{|\alpha|}$$ the canonical projection given by
$\pi_{\alpha}(\mathrm{I}) = (i_k)_{k \in \alpha}$ for any
$\mathrm{I} = (i_1, \ldots, i_d) \in [n]^d$. This notation allows
us express the inequality above in a much more understandable way.
Namely, we have

\begin{equation} \label{Equation-Iteration}
\|\mathcal{S}_d(a)\|_{L_p(\varphi \otimes \tau_d)} \le 2^d
\max_{(\alpha,\beta) \in \mathbb{P}_d(2)} \left\{ \Big\|
\sum_{\mathrm{I} \in [n]^d} a_{\mathrm{I}}^{} \otimes
e_{\pi_{\alpha}(\mathrm{I}),\pi_{\beta}(\mathrm{I})}
\Big\|_{L_p(\varphi; S_p^{n^d})} \ \right\}.
\end{equation}

\begin{remark} \label{Remark-Converse}
\textnormal{By the same arguments, the converse of
(\ref{Equation-Iteration}) holds with constant 1.}
\end{remark}

\subsection{Fell's absorption principle in $L_p$}

In the following, we shall use repeatedly the following
$L_p$-valued version of the so-called Fell's absorption principle.
This result might be known as folklore in the theory. However, we
include the proof since we were not able to provide a reference.

\vskip5pt

\noindent \textbf{Absorption Principle in $L_p$.} \emph{Given a
discrete group $\mathrm{G}$, let us denote by
$\lambda_{\mathrm{G}}$ the left regular representation of
$\mathrm{G}$ and by $\tau_{\mathrm{G}}$ the associated trace on
the reduced $C^*$-algebra of $\mathrm{G}$. Then, given any other
unitary representation $\pi: \mathrm{G} \rightarrow
\pi(\mathrm{G})''$, the following representations are unitarily
equivalent $$\lambda_{\mathrm{G}} \otimes \pi \simeq
\lambda_{\mathrm{G}} \otimes 1,$$ where $1$ stands for the trivial
representation of $\mathrm{G}$ in $\pi(\mathrm{G})''$. Moreover,
let us take any faithful normalized trace $\psi$ on
$\pi(\mathrm{G})''$. Then, given any function $a: \mathrm{G}
\rightarrow L_p(\varphi)$ finitely supported on $\mathrm{G}$, the
following equality holds for $1 \le p \le \infty$}
\begin{equation} \label{Equation-Fell}
\Big\| \sum_{t \in \mathrm{G}} a(t) \otimes
\lambda_{\mathrm{G}}(t) \otimes \pi(t) \Big\|_{L_p(\varphi \otimes
\tau_{\mathrm{G}} \otimes \psi)} = \Big\| \sum_{t \in \mathrm{G}}
a(t) \otimes \lambda_{\mathrm{G}}(t) \Big\|_{L_p(\varphi \otimes
\tau_{\mathrm{G}})}.
\end{equation}

\dem We refer the reader to Proposition 8.1 of \cite{P2} for a
proof of the claimed unitary equivalence. For the second
assertion, it is easy to reduce to the case when $\varphi$ is a
tracial state. Then, we fix a pair of operators
\begin{eqnarray*}
\mathsf{S} & = & \sum_{t \in \mathrm{G}} a(t) \otimes
\lambda_{\mathrm{G}}(t), \\ \mathsf{T} & = & \sum_{t \in
\mathrm{G}} a(t) \otimes \lambda_{\mathrm{G}}(t) \otimes \pi(t),
\end{eqnarray*}
with $a: \mathrm{G} \rightarrow L_p(\varphi) \cap
L_\infty(\varphi)$ finitely supported. Since $\psi$ is normalized,
it is clear that
$$(\varphi \otimes \tau_{\mathrm{G}}) (\mathsf{S}) = \varphi \big(
a(e) \big) = \varphi \big( a(e) \big) \psi \big( \pi(e) \big) =
(\varphi \otimes \tau_{\mathrm{G}} \otimes \psi) (\mathsf{T}).$$
Let us consider the operators $x = \mathsf{S}^* \mathsf{S}$ and $z
= \mathsf{T}^* \mathsf{T}$. Recalling that the equality above
holds for any pair of operators of the same kind, we deduce that
$$(\varphi \otimes \tau_{\mathrm{G}})(x^n) = (\varphi \otimes
\tau_{\mathrm{G}} \otimes \psi)(z^n)$$ for any integer $n \ge 0$.
Now we let $\mathrm{A}_x$ (resp. $\mathrm{A}_z$) be the
(commutative) algebra generated by $x$ (resp. $z$) in
$L_{\infty}(\varphi \otimes \tau_{\mathrm{G}})$ (resp.
$L_{\infty}(\varphi \otimes \tau_{\mathrm{G}} \otimes \psi)$). If
$\mu_x$ (resp. $\mu_z$) denotes the inherited probability measure
on $\mathrm{A}_x$ (resp. $\mathrm{A}_z$), we have $$\int
\mathrm{Q}(x) \, d \mu_x = \int \mathrm{Q}(z) \, d \mu_z,$$ for
any polynomial $\mathrm{Q}$. By the Stone-Weierstrass theorem, we
conclude that the distribution of $x$ with respect to $\mu_x$
coincides with the distribution of $z$ with respect to $\mu_z$.
Therefore, $\|\mathsf{S}\|_{L_p(\varphi \otimes
\tau_{\mathrm{G}})} = \|x\|_{L_{p/2}(\mu_x)} =
\|z\|_{L_{p/2}(\mu_z)} = \|\mathsf{T}\|_{L_p(\varphi \otimes
\tau_{\mathrm{G}} \otimes \psi)}$. \fin

\section{Khintchine type inequalities for $\Word_p(n,d)$}
\label{Section-Main}

We first prove the case of degree $2$, since the notation is
simpler and it contains almost all the ingredients employed in the
proof of the general case. This will simplify the reading of the
paper. After the proof of the general case, we study the same
problem when redefining the spaces $\Max_p(n,d)$ and
$\Word_p(n,d)$ so that we consider all the reduced words of length
$d$.

\subsection{The case of degree $2$}
\label{Subsection-d=2}

As we did in the Introduction, we define $\Word_p(n,2)$ to be the
subspace of $L_p(\tau)$ generated by the operators
$\lambda(g_ig_j)$ for $1 \le i,j \le n$. We shall also consider
the spaces $$\begin{array}{lcl} \displaystyle \Max_p(n,2) =
\sum_{k=0}^2 C_p^{n^k} \otimes_h R_p^{n^{d-k}} & \qquad &
\mbox{for} \ 1 \le p \le 2, \\ \displaystyle \Max_p(n,2) =
\bigcap_{k=0}^2 C_p^{n^k} \otimes_h R_p^{n^{d-k}} & \qquad &
\mbox{for} \ 2 \le p \le \infty.
\end{array}$$ That is, if $\mathcal{A} = \big\{ a_{ij} \, | \
1 \le i,j \le n \big\} \subset L_p(\varphi)$, we consider the
norms $$\Big\| \sum_{i,j=1}^n a_{ij} \otimes e_{1,ij}
\Big\|_{L_p(\varphi; R_p^{n^2})}, \quad \Big\| \sum_{i,j=1}^n
a_{ij} \otimes e_{ij} \Big\|_{L_p(\varphi; S_p^n)}, \quad \Big\|
\sum_{i,j=1}^n a_{ij} \otimes e_{ij,1} \Big\|_{L_p(\varphi;
C_p^{n^2})}.$$ We label them by $\|\mathcal{A}\|_0,
\|\mathcal{A}\|_1$ and $\|\mathcal{A}\|_2$ respectively. Then, we
have $$\|\mathcal{A}\|_{L_p(\varphi; \Max_p(n,2))} = \, \max
\left\{ \|\mathcal{A}\|_k \, \Big| \ 0 \le k \le 2 \right\} \qquad
\mbox{for} \ 2 \le p \le \infty.$$ This identity describes the
operator space structure of $\Max_p(n,2)$ for $2 \le p \le
\infty$. On the other hand, the obvious modifications lead to a
description of the operator space structure of $\Max_p(n,2)$ for
$1 \le p \le 2$. We shall prove the following result

\begin{theorem} \label{Theorem-Degree2}
$\Max_p(n,2)$ and $\Word_p(n,2)$ are completely isomorphic for $1
\le p \le \infty$. More concretely, there exists an absolute
constant $\mathrm{c}$ independent of $n$ and $p$ such that the
following inequalities hold $$\frac{1}{\mathrm{c}} \,
\|\mathcal{A}\|_{L_p(\varphi; \Max_p(n,2))} \le \Big\|
\sum_{i,j=1}^n a_{ij} \otimes \lambda(g_ig_j) \Big\|_{L_p(\varphi
\otimes \tau)} \le \mathrm{c} \, \|\mathcal{A}\|_{L_p(\varphi;
\Max_p(n,2))}.$$ Moreover, the natural projection $\mathrm{P}:
L_p(\tau) \rightarrow \Word_p(n,2)$ is c.b. with
$\|\mathrm{P}\|_{cb} \le \mathrm{c}$.
\end{theorem}

A similar statement to this result was proved by Haagerup
(unpublished) and Buchholz \cite{Bu1} for $p=\infty$, see also
Theorem 9.7.4 in \cite{P2}. Namely, the only difference is that
Buchholz considered the whole set of words of length $2$ instead
of those words composed only of generators. Therefore, it is clear
that Theorem \ref{Theorem-Degree2} holds for $p=\infty$. By
transposition, $\mathrm{P}$ also defines a completely bounded
projection from $L_1(\tau)$ onto $\Word_1(n,2)$. Hence, the last
assertion of Theorem \ref{Theorem-Degree2} follows by complex
interpolation. In particular, if $p'$ stands for the conjugate
exponent of $p$, it turns out that the dual $\Word_p(n,2)^*$ is
completely isomorphic to $\Word_{p'}(n,2)$ for $1 \le p \le
\infty$. On the other hand, it is obvious that $\Max_p(n,2)^*$ is
completely isometric to $\Max_{p'}(n,2)$. In summary, it suffices
to prove Theorem \ref{Theorem-Degree2} for $2 \le p \le \infty$
since the case $1 \le p \le 2$ follows by duality.

\begin{remark}
\textnormal{From the previous considerations, it is clear that the
space $\Word_p(n,2)$ is completely isomorphic to
$(\Word_{\infty}(n,2), \Word_2(n,2))_{\theta}$ for $\theta = 2/p$.
Hence, if we knew a priori that $\Max_p(n,2)$ is an interpolation
family for $2 \le p \le \infty$, Theorem \ref{Theorem-Degree2}
would follow by complex interpolation between the obvious case for
$p=2$ and Buchholz's result. Conversely, Theorem
\ref{Theorem-Degree2} implies that $\Max_p(n,2)$ is an
interpolation family for $2 \le p \le \infty$. The fact that
$\Max_p(n,2)$ is an interpolation family was known to M. Junge
\cite{J1} at the time of the preparation of this paper. He
communicated to us the following sketch of his argument. Let
$\tau_n$ stand for the normalized trace on the matrix algebra
$M_n$ and let us write $\mathcal{N}$ for the free product of
algebras $$\mathcal{N} = \astop_{k=1}^n \mathcal{A}_k \qquad
\mbox{with} \qquad \mathcal{A}_k = (M_n \oplus_{\infty} M_n,
\varphi_n) \qquad \mbox{and} \qquad \varphi_n = \frac{1}{2}
(\tau_n \oplus \tau_n)$$ for $1 \le k \le n$. Then, if $\pi_k:
\mathcal{A}_k \rightarrow \mathcal{N}$ denotes the natural
inclusion, the map $$x \in R_p^{n^2} \cap S_p^n \cap C_p^{n^2}
\longmapsto \sum_{k=1}^n \pi_k((x,-x)) \in L_p(\mathcal{N})$$ is a
complete isomorphism onto its image and the image is completely
complemented in $L_p(\mathcal{N})$. Moreover, the constants
appearing in the complete isomorphism and the projection
considered above do not depend on $n$. This is based on the $L_p$
version of the operator-valued Voiculescu's inequality given in
\cite{J2}. Here we shall give a different proof that will be
useful in the proof of the general case of degree $d$.}
\end{remark}

Now we focus on the proof for the case $2 \le p \le \infty$. The
lower estimate is much simpler and it even holds with
$\mathrm{c}=1$. Namely, it suffices to check that
$$\|\mathcal{A}\|_k \le \Big\| \sum_{i,j=1}^n a_{ij} \otimes
\lambda(g_ig_j) \Big\|_{L_p(\varphi \otimes \tau)},$$ for any $0
\le k \le 2$. But we know that it holds trivially for $p=2$ and
also, by Buchholz's result, for $p=\infty$. Therefore, the lower
estimate follows by complex interpolation since $R_p^{n^2}$,
$S_p^n$ and $C_p^{n^2}$ are interpolation families. Hence, we just
need to prove the upper estimate. To that aim, we go back to
Section \ref{Section-Iteration}, where we considered the group
$\mathrm{G}_2 = \Free_n \times \Free_n$ and the subspace
$\Word_p(n) \otimes \Word_p(n)$ of $L_p(\tau_2)$ generated by the
family of operators $$\lambda(g_i) \otimes \lambda(g_j).$$ We
consider the subspace $\mathcal{V}_p(n,2)$ of $L_p(\tau_3)$
defined by $$\mathcal{V}_p(n,2) = \Big\{ \sum_{i,j=1}^n
\alpha_{ij} \, \lambda(g_i) \otimes \lambda(g_j) \otimes
\lambda(g_ig_j) \in L_p(\tau_3) \, \big| \ \alpha_{ij} \in
\mathbb{C} \Big\}.$$

\begin{lemma} \label{Lemma-Projection} $\mathcal{V}_p(n,2)$ is a
completely complemented subspace of $L_p(\tau_3)$.
\end{lemma}

\dem We know that both $\Word_p(n)$ and $\Word_p(n,2)$ are
completely complemented in $L_p(\tau)$. In particular, the
projection which maps $$\sum_{u,v,w \in \Free_n} \alpha_{uvw} \,
\lambda(u) \otimes \lambda(v) \otimes \lambda(w) \in L_p(\tau_3)$$
to the sum $$\Sigma_2(\alpha) = \sum_{i,j,r,s=1}^n \alpha_{ijrs}
\, \lambda(g_i) \otimes \lambda(g_j) \otimes \lambda(g_rg_s) \in
\Word_p(n) \otimes \Word_p(n) \otimes \Word_p(n,2),$$ is
completely bounded with cb norm uniformly bounded in $n$ and $p$.
This shows that $\Word_p(n) \otimes \Word_p(n) \otimes
\Word_p(n,2)$ is completely complemented in $L_p(\tau_3)$. After
that, we project onto $\mathcal{V}_p(n,2)$ by using the standard
diagonal projection $$\mathrm{P} \big( \Sigma_2(\alpha) \big) =
\sum_{i,j,r,s=1}^n \int \int \varepsilon_i \delta_j \Big[
\alpha_{ijrs} \, \lambda(g_i) \otimes \lambda(g_j) \otimes
\lambda(g_rg_s) \Big] \varepsilon_r \delta_s \, d\mu(\varepsilon)
d\mu(\delta),$$ where $\mu$ is the normalized counting measure on
$\{-1,1\}^n$. Now, in the case $p=\infty$ it is not difficult to
see that the norm of $$\Sigma_2^{\pm}(\alpha) = \sum_{i,j,r,s=1}^n
\varepsilon_i \delta_j \Big[ \alpha_{ijrs} \, \lambda(g_i) \otimes
\lambda(g_j) \otimes \lambda(g_rg_s) \Big] \varepsilon_r
\delta_s$$ in $L_{\infty}(\tau_3)$ is equivalent (in the category
of operator spaces) to that of $\Sigma_2(\alpha)$ in
$L_{\infty}(\tau_3)$ for any choice of signs $\varepsilon_i,
\delta_j, \varepsilon_r, \delta_s$. Moreover, the constants do not
depend on the signs taken. Indeed, by Buchholz's result the norm
of $\Sigma_2^{\pm}(\alpha)$ is equivalent to the norm of an element in
$L_{\infty}(\tau_2; \Max_{\infty}(n,d))$. In that case, the signs
$\varepsilon_r \delta_s$ can be regarded as a Schur multiplier.
Hence, $\varepsilon_r$ and $\delta_s$ can be dropped by means of
\cite[Exercise 1.5]{P2}. On the other hand, the signs
$\varepsilon_i$ and $\delta_j$ disappear by applying Fell's
absorption principle. Moreover, the norms of $\Sigma_2(\alpha)$
and $\Sigma_2^{\pm}(\alpha)$ clearly coincide for $p=2$. Therefore,
since $\Word_p(n) \otimes \Word_p(n) \otimes \Word_p(n,2)$ is an
interpolation family, both norms are equivalent for any $2 \le p
\le \infty$. In particular, by Jensen's inequality we have
$$\big\| \mathrm{P} \big( \Sigma_2(\alpha) \big) \big
\|_{L_p(\tau_3)} \le \mathrm{c} \, \big\| \Sigma_2(\alpha)
\big\|_{L_p(\tau_3)},$$ for some absolute constant $\mathrm{c}$.
Since the same holds taking values in $S_p$, it turns out that
$\mathrm{P}$ is a completely bounded projection with constants
independent of $n$ and $p$. In summary, putting all together the
result follows. This completes the proof. \fin

The next step in the proof of Theorem \ref{Theorem-Degree2} is to
show that $\Word_p(n,2)$ and $\mathcal{V}_p(n,2)$ are completely
isomorphic operator spaces for any exponent $2 \le p \le \infty$.
Namely, given a family $\mathcal{A} = \big\{ a_{ij} \, | \ 1 \le
i,j \le n \big\}$ in $L_p(\varphi)$, we consider the sum
$$\Sigma_2(a) = \sum_{i,j=1}^n \lambda(g_i) \otimes \lambda(g_j)
\otimes a_{ij} \otimes \lambda(g_ig_j) \in L_p(\tau_2 \otimes
\varphi \otimes \tau).$$ Applying Buchholz's result to
$\Sigma_2(a)$, we get the equivalence of norms $$\big\|
\Sigma_2(a) \big\|_{L_{\infty}(\tau_2 \otimes \varphi \otimes
\tau)} \simeq \max \left\{ \|\mathcal{A}'\|_k \, \Big| \ 0 \le k
\le 2 \right\},$$ where $\mathcal{A}' = \big\{ a_{ij}' \, | \ 1
\le i,j \le n \big\}$ with $$a_{ij}' = \lambda(g_i) \otimes
\lambda(g_j) \otimes a_{ij}.$$

\begin{remark} \label{Remark-Factorization-Unitary-1}
\textnormal{Note that $$\sum_{i,j=1}^n a_{ij}' \otimes e_{ij} =
\Phi_1 \cdot \Big[ \sum_{i,j=1}^n 1 \otimes 1 \otimes a_{ij}
\otimes e_{ij} \Big] \cdot \Phi_2,$$ where
\begin{eqnarray*}
\Phi_1 & = & \sum_{i=1}^n \lambda(g_i) \otimes 1 \otimes 1 \otimes
e_{ii}, \\ \Phi_2 & = & \sum_{j=1}^n 1 \otimes \lambda(g_j)
\otimes 1 \otimes e_{jj}.
\end{eqnarray*}}
\end{remark}

Therefore, according to Remark
\ref{Remark-Factorization-Unitary-1} and since $\Phi_1$ and
$\Phi_2$ are unitary, we obtain $\|\mathcal{A}\|_1 =
\|\mathcal{A}'\|_1$. The obvious modifications lead to
$\|\mathcal{A}\|_k = \|\mathcal{A}'\|_k$ for $k=0$ and $k=2$. In
summary, if $$\mathcal{S}_2(a) = \sum_{i,j=1}^n a_{ij} \otimes
\lambda(g_ig_j),$$ we conclude that the norm of $\mathcal{S}_2(a)$
in $L_{\infty}(\varphi \otimes \tau)$ is equivalent (in the
category of operator spaces) to the norm of $\Sigma_2(a)$ in
$L_{\infty}(\tau_2 \otimes \varphi \otimes \tau)$. Now recall that
by Lemma \ref{Lemma-Projection}, $\mathcal{V}_p(n,2)$ is an
interpolation family for $2 \le p \le \infty$. Then, since the
norm of these sums obviously coincide when $p=2$, we get by
complex interpolation
\begin{equation} \label{Equation-Tensors}
\Big\| \sum_{i,j=1}^n a_{ij} \otimes \lambda(g_ig_j) \Big\|_p \le
\mathrm{c} \Big\| \sum_{i,j=1}^n \lambda(g_i) \otimes \lambda(g_j)
\otimes a_{ij} \otimes \lambda(g_ig_j) \Big\|_p
\end{equation}
for any $2 \le p \le \infty$. Here $\mathrm{c}$ denotes an
absolute constant independent of $n$ and $p$. In what follows, the
value of $\mathrm{c}$ might change from one instance to another.
Now we apply the iteration of the Khintchine inequality
(\ref{Equation-Iteration}) to inequality (\ref{Equation-Tensors})
to obtain
\begin{equation} \label{Equation-2Partition}
\Big\| \sum_{i,j=1}^n a_{ij} \otimes \lambda(g_ig_j) \Big\|_p \le
\mathrm{c} \, \max_{(\alpha,\beta) \in \mathbb{P}_2(2)} \left\{
\Big\| \sum_{\mathrm{I} \in [n]^2} \tilde{a}_{\mathrm{I}}^{}
\otimes e_{\pi_{\alpha}(\mathrm{I}), \pi_{\beta}(\mathrm{I})}
\Big\|_{L_p(\varphi \otimes \tau; S_p^{n^2})} \right\},
\end{equation}
with $\tilde{a}_{ij} = a_{ij} \otimes \lambda(g_ig_j)$. Hence, we
have four terms on the right $$\begin{array}{lcl} \bullet \
\displaystyle \Big\| \sum_{i,j=1}^n \tilde{a}_{ij} \otimes
e_{1,ij} \Big\|_{L_p(\varphi \otimes \tau; R_p^{n^2})} & \qquad &
\bullet \ \displaystyle \Big\| \sum_{i,j=1}^n \tilde{a}_{ij}
\otimes e_{ij} \Big\|_{L_p(\varphi \otimes \tau; S_p^n)} \\
\bullet \ \displaystyle \Big\| \sum_{i,j=1}^n \tilde{a}_{ij}
\otimes e_{ij,1} \Big\|_{L_p(\varphi \otimes \tau; C_p^{n^2})} &
\qquad & \bullet \ \displaystyle \Big\| \sum_{i,j=1}^n
\tilde{a}_{ij} \otimes e_{ji} \Big\|_{L_p(\varphi \otimes \tau;
S_p^n)} \end{array}$$ If $\widetilde{\mathcal{A}} = \{
\tilde{a}_{ij} \, | \ 1 \le i,j \le n \}$, the first three terms
are nothing but $\|\widetilde{\mathcal{A}}\|_0,
\|\widetilde{\mathcal{A}}\|_1, \|\widetilde{\mathcal{A}}\|_2$.
Arguing as above we have $$\max \left\{ \|\mathcal{A}\|_k \, \Big|
\ 0 \le k \le 2 \right\} = \max \left\{
\|\widetilde{\mathcal{A}}\|_k \, \Big| \ 0 \le k \le 2 \right\}.$$
In particular, the proof will be completed if we see that $$\Big\|
\sum_{i,j=1}^n \tilde{a}_{ij} \otimes e_{ji} \Big\|_{L_p(\varphi
\otimes \tau; S_p^n)} \le \mathrm{c} \, \max \left\{
\|\mathcal{A}\|_k \, \Big| \ 0 \le k \le 2 \right\}.$$ This is the
content of the following Lemma. The proof is not complicated but,
as we shall see in the next paragraph, it constitutes one of the
key points in the proof of the general case.

\begin{lemma} \label{Lemma-Transposition}
The following inequality holds
\begin{eqnarray*}
\lefteqn{\Big\| \sum_{i,j=1}^n a_{ij} \otimes \lambda(g_ig_j)
\otimes e_{ji} \Big\|_{L_p(\varphi \otimes \tau; S_p^n)}} \\ & \le
& \mathrm{c} \, \max \left\{ \Big\| \sum_{i,j=1}^n a_{ij} \otimes
e_{1,ij} \Big\|_{L_p(\varphi;R_p^{n^2})} , \Big\| \sum_{i,j=1}^n
a_{ij} \otimes e_{ij,1} \Big\|_{L_p(\varphi;C_p^{n^2})} \right\}.
\end{eqnarray*}
\end{lemma}

\dem When $p=\infty$ we can apply Buchholz's result to obtain
$$\Big\| \sum_{i,j=1}^n a_{ij} \otimes \lambda(g_ig_j) \otimes
e_{ji} \Big\|_{L_{\infty}(\varphi \otimes \tau; S_{\infty}^n)} \le
\mathrm{c} \, \max \Big\{ \mathrm{A}, \mathrm{B}, \mathrm{C}
\Big\},$$ where the terms $\mathrm{A}$ and $\mathrm{C}$ are given
by
\begin{eqnarray*}
\mathrm{A} & = & \Big\| \Big( \sum_{i,j=1}^n (a_{ij} \otimes
e_{ji}) (a_{ij} \otimes e_{ji})^* \Big)^{1/2} \Big\|_{\infty} =
\sup_{1 \le j \le n} \Big\| \Big( \sum_{i=1}^n a_{ij} a_{ij}^*
\Big)^{1/2} \Big\|_{\infty} \\ \mathrm{C} & = & \Big\| \Big(
\sum_{i,j=1}^n (a_{ij} \otimes e_{ji})^* (a_{ij} \otimes e_{ji})
\Big)^{1/2} \Big\|_{\infty} = \sup_{1 \le i \le n} \Big\| \Big(
\sum_{j=1}^n a_{ij}^* a_{ij} \Big)^{1/2} \Big\|_{\infty}.
\end{eqnarray*}
In particular, we have the following estimates
\begin{eqnarray*}
\mathrm{A} & \le & \Big\| \Big( \sum_{i,j=1}^n a_{ij} a_{ij}^*
\Big)^{1/2} \Big\|_{\infty} = \Big\| \sum_{i,j=1}^n a_{ij} \otimes
e_{1,ij} \Big\|_{L_{\infty}(\varphi; R_{\infty}^{n^2})} \\
\mathrm{C} & \le & \Big\| \Big( \sum_{i,j=1}^n a_{ij}^* a_{ij}
\Big)^{1/2} \Big\|_{\infty} = \Big\| \sum_{i,j=1}^n a_{ij} \otimes
e_{ij,1} \Big\|_{L_{\infty}(\varphi; C_{\infty}^{n^2})}.
\end{eqnarray*}
It remains to estimate the middle term $\mathrm{B}$. We have
\begin{eqnarray*}
\mathrm{B} & = & \Big\| \sum_{i,j=1}^n a_{ij} \otimes e_{ji}
\otimes e_{ij} \Big\|_{L_{\infty}(\varphi; S_{\infty}^n
\otimes_{\min} S_{\infty}^n)} \\ & = & \Big\| \Big( \sum_{i,j=1}^n
a_{ij}^* a_{ij} \otimes e_{ii} \otimes e_{jj} \Big)^{1/2}
\Big\|_{L_{\infty}(\varphi; S_{\infty}^n \otimes_{\min}
S_{\infty}^n)} \\ & = & \sup_{1 \le i,j \le n}
\|a_{ij}\|_{L_{\infty}(\varphi)}.
\end{eqnarray*}
Therefore, the term $\mathrm{B}$ is even smaller that
$\mathrm{A}$, $\mathrm{C}$. This completes the proof for the case
$p=\infty$. On the other hand, for the case $p=2$ we clearly have
an equality. Finally, we recall that $$L_p \big( \varphi;
S_p^n(\Word_p(n,2))^{op} \big) \qquad \mbox{and} \qquad L_p
\big(\varphi; R_p^{n^2} \cap C_p^{n^2} \big)$$ are interpolation
families, see Section 9.5 in \cite{P2} for the details. Therefore,
the result follows for $2 \le p \le \infty$ by complex
interpolation. This completes the proof. \fin

\subsection{The general case}
\label{Subsection-General}

Now we prove the analog of Theorem \ref{Theorem-Degree2} for any
positive integer $d$. As we shall see, there exist a lot of
similarities with the proof for degree $2$. Therefore, we shall
not repeat in detail those arguments which already appeared above.
The statement of this result is the following.

\begin{theorem} \label{Theorem-Degreed}
$\Max_p(n,d)$ and $\Word_p(n,d)$ are completely isomorphic for $1
\le p \le \infty$. More concretely, there exists an absolute
constant $\mathrm{c}_d$ depending only on $d$ such that the
following inequalities hold $$\frac{1}{\mathrm{c}_d} \,
\|\mathcal{A}\|_{L_p(\varphi; \Max_p(n,d))} \le \Big\| \sum_{i_1
\ldots i_d=1}^n a_{i_1 \cdots i_d} \otimes \lambda(g_{i_1} \cdots
g_{i_d}) \Big\|_{L_p(\varphi \otimes \tau)} \le \mathrm{c}_d \,
\|\mathcal{A}\|_{L_p(\varphi; \Max_p(n,d))}.$$ Moreover, the
natural projection $\mathrm{P}: L_p(\tau) \rightarrow
\Word_p(n,d)$ is c.b. with $\|\mathrm{P}\|_{cb} \le \mathrm{c}_d$.
\end{theorem}

Before starting the proof of Theorem \ref{Theorem-Degreed}, we
point out some remarks analogous to those given for Theorem
\ref{Theorem-Degree2}. The arguments needed to prove the
assertions given below are the same as the ones we used for the
case of degree $2$.
\begin{itemize}
\item Again, a similar statement to this result was
proved by Buchholz \cite{Bu1} for $p=\infty$. Buchholz's
considered the whole set of words of length $d$. Therefore,
Theorem \ref{Theorem-Degreed} holds for $p=\infty$ by Buchholz's
more general statement.

\item By transposition and complex
interpolation, the last assertion of Theorem \ref{Theorem-Degreed}
follows. Hence, $\Word_p(n,d)$ interpolates well up to complete
isomorphism for $1 \le p \le \infty$. Moreover, $\Word_{p'}(n,d)$
is completely isomorphic to $\Word_p(n,d)^*$. In particular, since
$\Max_p(n,d)$ behaves well with respect to duality, it suffices to
prove Theorem \ref{Theorem-Degreed} for $2 \le p \le \infty$.

\item Given a family of operators $\mathcal{A} = \big\{
a_{i_1 \cdots i_d} \, | \ 1 \le i_k \le n \big\}$ in
$L_p(\varphi)$, we define $$\|\mathcal{A}\|_k = \Big\| \sum_{i_1,
\ldots, i_d=1}^n a_{i_1 \cdots i_d} \otimes e_{(i_1 \cdots i_k),
(i_{k+1} \cdots i_d)} \Big\|_{L_p \left( \varphi; C_p^{n^k}
\otimes_h R_p^{n^{d-k}} \right)}.$$

\item If $2 \le p \le \infty$, the lower estimate
holds with $\mathrm{c}_d = 1$. Namely, it follows by Buchholz's
result and complex interpolation since we are allowed to look
separately at the inequalities $$\qquad \|\mathcal{A}\|_k \le
\Big\| \sum_{i_1, \ldots, i_d=1}^n a_{i_1 \cdots i_d} \otimes
\lambda(g_{i_1} \cdots g_{i_d}) \Big\|_{L_p(\varphi \otimes
\tau)}, \qquad \mbox{for} \ 0 \le k \le d.$$ Recall that
$C_p^{n^k} \otimes_h R_p^{n^{d-k}}$ is a rectangular Schatten
$p$-class of size $n^k \times n^{d-k}$. In particular, it follows
that it is an interpolation family.
\end{itemize}

In what follows, $\mathrm{c}_d$ will denote a constant depending
only on $d$ and whose value might change from one instance to
another. Now we start the proof of Theorem \ref{Theorem-Degreed}.
By the considerations above, it suffices to prove the upper
estimate for $2 \le p \le \infty$. Let $\mathrm{G}$ stand for the
free group $\Free_{nd}$ with $nd$ generators and let $\psi_{nd}$
be the natural trace on the reduced $C^*$-algebra of $\mathrm{G}$.
We label the generators by $g_{1k}, g_{2k}, \ldots, g_{dk}$ with
$1 \le k \le n$. If $\lambda_{\mathrm{G}}$ denotes the left
regular representation of $\mathrm{G}$, we consider the family of
operators $$\mathcal{A}^{\dag} = \Big\{ a_{i_1i_2 \cdots i_d}
\otimes \lambda_{\mathrm{G}}(g_{1i_1} g_{2i_2} \cdots g_{di_d}) \,
\Big| \ 1 \le i_k \le n \Big\}.$$ That is, we take the image under
$\lambda_{\mathrm{G}}$ of the set of reduced words of length $d$
where the first letter is one of the first $n$ generators, the
second letter is one of the second $n$ generators and so on. Let
us write $\Word_p(nd,n,d)$ to denote the subspace of
$L_p(\psi_{nd})$ generated by the operators $\lambda_{\mathrm{G}}
(g_{1i_1} \cdots g_{di_d})$. Then, we have
\begin{equation} \label{Equation-Independence}
\Big\| \sum_{i_1, \ldots, i_d=1}^n a_{i_1 \cdots i_d} \otimes
\lambda(g_{i_1} \cdots g_{i_d}) \Big\|_p \le \mathrm{c}_d  \Big\|
\sum_{i_1, \ldots, i_d=1}^n a_{i_1 \cdots i_d} \otimes
\lambda_{\mathrm{G}}(g_{1i_1} \cdots g_{di_d}) \Big\|_p.
\end{equation}
Namely, for $p=\infty$ this follows by Buchholz's result. Then,
since both $\Word_p(n,d)$ and $\Word_p(nd,n,d)$ are interpolation
families, inequality (\ref{Equation-Independence}) holds for any
$2 \le p \le \infty$ by complex interpolation. Now, proceeding as
in the previous paragraph, we consider the group $\mathrm{G}_d =
\Free_n \times \Free_n \times \cdots \times \Free_n$ and the
subspace $\Word_p(n)^{\otimes d}$ of $L_p(\tau_d)$ generated by
the family of operators $$\lambda(g_{i_1}) \otimes
\lambda(g_{i_2}) \otimes \cdots \otimes \lambda(g_{i_d}).$$ Then
we define $\mathcal{V}_p(n,d)$ as the subspace of $L_p(\tau_d
\otimes \psi_{nd})$ defined by $$\mathcal{V}_p(n,d) = \Big\{
\sum_{i_1, \ldots, i_d=1}^n \alpha_{i_1 \cdots i_d} \,
\lambda(g_{i_1}) \otimes \cdots \otimes \lambda(g_{i_d}) \otimes
\lambda_{\mathrm{G}}(g_{1i_1} \cdots g_{di_d}) \, \big| \
\alpha_{i_1 \cdots i_d} \in \mathbb{C} \Big\}.$$ Recalling that
both $\Word_p(n)$ and $\Word_p(nd,n,d)$ are completely
complemented in their respective $L_p$ spaces, it can be showed
just like in Lemma \ref{Lemma-Projection} that
$\mathcal{V}_p(n,d)$ is completely complemented in $L_p(\tau_d
\otimes \psi_{nd})$ with constants depending only on the degree
$d$. The next step in the proof is to obtain the analog of
inequality (\ref{Equation-Tensors}). Namely, the inequality
\begin{eqnarray} \label{Equation-d-Tensors}
\lefteqn{\Big\| \sum_{i_1, \ldots, i_d=1}^n a_{i_1 \cdots i_d}
\otimes \lambda_{\mathrm{G}}(g_{1i_1} \cdots g_{di_d}) \Big\|_p}
\\ \nonumber & \le & \mathrm{c}_d \, \Big\| \sum_{i_1, \ldots,
i_d=1}^n \lambda(g_{i_1}) \otimes \cdots \otimes \lambda(g_{i_d})
\otimes a_{i_1 \cdots i_d} \otimes \lambda_{\mathrm{G}}(g_{1i_1}
\cdots g_{di_d}) \Big\|_p,
\end{eqnarray}
for $2 \le p \le \infty$. The proof of this inequality is
identical to the one given for inequality
(\ref{Equation-Tensors}). Indeed, we have just showed that both
$\Word_p(nd,n,d)$ and $\mathcal{V}_p(n,d)$ are interpolation
families. Therefore, the proof of (\ref{Equation-d-Tensors}) works
by complex interpolation between the obvious case $p=2$ and the
case $p=\infty$. When $p=\infty$, the idea consists in applying to
both terms in (\ref{Equation-d-Tensors}) Buchholz's result for
degree $d$. Then, inequality (\ref{Equation-d-Tensors}) becomes
equivalent to $$\max \left\{ \|\mathcal{A}\|_k \, \Big| \ 0 \le k
\le d \right\} \le \mathrm{c}_d \, \max \left\{ \|\mathcal{A}'\|_k
\, \Big| \ 0 \le k \le d \right\},$$ where $\mathcal{A}'$ is given
by $$\mathcal{A}' = \Big\{ a_{i_1i_2 \cdots i_d} \otimes
\lambda(g_{i_1}) \otimes \cdots \otimes \lambda(g_{i_d}) \, \Big|
\ 1 \le i_k \le n \Big\}.$$

\begin{remark} \label{Remark-Factorization-Unitary-2}
\textnormal{Again it is clear that $\|\mathcal{A}\|_k =
\|\mathcal{A}'\|_k$ for any $0 \le k \le d$. Namely, as we pointed
out in Remark \ref{Remark-Factorization-Unitary-1}, the sum
$$\sum_{i_1, \ldots, i_d =1}^n \lambda(g_{i_1}) \otimes \cdots
\otimes \lambda(g_{i_d}) \otimes a_{i_1i_2 \cdots i_d} \otimes
e_{(i_1 \cdots i_k),(i_{k+1} \cdots i_d)}$$ factorizes as $$\Phi_1
\cdot \Big[ \sum_{i_1, \ldots, i_d =1}^n a_{i_1i_2 \cdots i_d}
\otimes e_{(i_1 \cdots i_k),(i_{k+1} \cdots i_d)} \Big] \cdot
\Phi_2,$$ where $\Phi_1$ (resp. $\Phi_2$) is a $n^k \times n^k$
(resp. $n^{d-k} \times n^{d-k}$) unitary mapping.}
\end{remark}

This completes the proof of inequality (\ref{Equation-d-Tensors}).
Another possible approach to (\ref{Equation-d-Tensors}) is given
by iterating Fell's absorption principle $d$ times, with the
suitable choice for $\pi$ each time. We leave the details to the
reader. Notice that Fell's absorption principle shows that
(\ref{Equation-d-Tensors}) is in fact an equality with
$\mathrm{c}_d = 1$. Then we apply the iteration of Khintchine
inequality (\ref{Equation-Iteration}) to inequality
(\ref{Equation-d-Tensors}). This gives
\begin{equation} \label{Equation-d-2Partition}
\Big\| \sum_{i_1, \ldots, i_d=1}^n a_{i_1 \cdots i_d}^{\dag}
\Big\|_p \le \mathrm{c}_d \, \max_{(\alpha,\beta) \in
\mathbb{P}_d(2)} \left\{ \Big\| \sum_{\mathrm{I} \in [n]^d}
a_{\mathrm{I}}^{\dag} \otimes e_{\pi_{\alpha}(\mathrm{I}),
\pi_{\beta}(\mathrm{I})} \Big\|_{L_p(\varphi \otimes \psi_{nd};
S_p^{n^d})} \right\},
\end{equation}
with $a_{i_1 \cdots i_d}^{\dag} = a_{i_1 \cdots i_d} \otimes
\lambda_{\mathrm{G}}(g_{1i_1} \cdots g_{di_d})$. Now we have $2^d$
terms on the right. Before going on, let us look for a moment at
the norm of the space $\Max_p(n,d)$. Concretely, if we rewrite the
definition of $\Max_p(n,d)$ for $2 \le p \le \infty$ with the
notation employed at the end of Section \ref{Section-Iteration},
we obtain
\begin{equation} \label{Equation-Kp(n,d)}
\|\mathcal{A}\|_{L_p(\varphi; \Max_p(n,d))} = \max_{0 \le k \le d}
\left\{ \Big\| \sum_{\mathrm{I} \in [n]^k} \sum_{\mathrm{J}
\in[n]^{d-k}} a_{\mathrm{IJ}}^{} \otimes e_{\mathrm{I},
\mathrm{J}}^{} \Big\|_{L_p \left( \varphi; C_p^{n^k} \otimes_h
R_p^{n^{d-k}} \right)} \right\}.
\end{equation}

To complete the proof of Theorem \ref{Theorem-Degreed} it remains
to see that the right side of (\ref{Equation-d-2Partition}) is
controlled by the right side of (\ref{Equation-Kp(n,d)}). Here is
where the proof of the general case differs from that of degree
$2$. Let us sketch briefly how we shall conclude the proof. If
$\mathrm{RHS}_{(\ref{Equation-d-2Partition})}$ stands for the
right hand side of (\ref{Equation-d-2Partition}), we shall prove
that
\begin{equation} \label{Equation-Plan}
\mathrm{RHS}_{(\ref{Equation-d-2Partition})} \le \mathrm{c}_d \
\max_{0 \le k \le d} \left\{ \Big\| \sum_{\mathrm{I} \in [n]^k}
\sum_{\mathrm{J} \in[n]^{d-k}} a_{\mathrm{IJ}}^{\dag} \otimes
e_{\mathrm{I}, \mathrm{J}}^{} \Big\|_{L_p \left( \varphi \otimes
\psi_{nd}; C_p^{n^k} \otimes_h R_p^{n^{d-k}} \right)} \right\}.
\end{equation}

Assuming we have (\ref{Equation-Plan}), the proof is completed
since the right side of (\ref{Equation-Plan}) coincides with the
right side of (\ref{Equation-Kp(n,d)}). Namely, it clearly follows
by the same factorization argument as above. That is,
$$\sum_{\mathrm{I} \in [n]^k} \sum_{\mathrm{J} \in[n]^{d-k}}
a_{\mathrm{IJ}}^{\dag} \otimes e_{\mathrm{I}, \mathrm{J}}^{} =
\Phi_1 \cdot \Big[ \sum_{\mathrm{I} \in [n]^k} \sum_{\mathrm{J}
\in[n]^{d-k}} a_{\mathrm{IJ}}^{} \otimes e_{\mathrm{I},
\mathrm{J}}^{} \Big] \cdot \Phi_2,$$ with $\Phi_1$ and $\Phi_2$
suitably chosen unitary mappings. In summary, it remains to prove
inequality (\ref{Equation-Plan}). Recall that the $d+1$ terms
which appear on the right hand side of (\ref{Equation-Plan}) also
appear on its left hand side. They correspond to $\alpha =
\emptyset, [1], [2], \ldots, [d]$. The remaining terms correspond
to certain transpositions just like the term we bounded with the
aid of Lemma \ref{Lemma-Transposition}. Thus, we just need to show
that the transposed terms are controlled by the non-transposed
ones. To that aim, we need to introduce some notation. Given
$(\alpha, \beta) \in \mathbb{P}_d(2)$, we define
\begin{eqnarray*}
\mathrm{a} & = & \max \Big\{ k \, | \ k \in \alpha \Big\}, \\
\mathrm{b} & = & \min \, \Big\{ k \, | \ k \in \beta \Big\}.
\end{eqnarray*}
We define $\mathrm{a} = 0$ for $\alpha = \emptyset$ and
$\mathrm{b} = d+1$ for $\beta = \emptyset$. We shall say that
$(\alpha, \beta)$ is \emph{non-transposed} whenever $\mathrm{a} <
\mathrm{b}$ and $(\alpha, \beta)$ will be called \emph{transposed}
otherwise. Let us introduce the number $\mathrm{T}(\alpha,\beta) =
\mathrm{a} - \mathrm{b}$, so that $(\alpha, \beta)$ is transposed
whenever $\mathrm{T}(\alpha, \beta) > 0$. The proof of the
remaining inequality lies on the following claim.

\begin{claim} \label{Claim}
Let $(\alpha, \beta)$ be a transposed element of
$\mathbb{P}_d(2)$. Then, we have
\begin{eqnarray*}
\lefteqn{\Big\| \sum_{\mathrm{I} \in [n]^d} a_{\mathrm{I}}^{\dag}
\otimes e_{\pi_{\alpha}(\mathrm{I}), \pi_{\beta}(\mathrm{I})}
\Big\|_p}  \\ & \le & \mathrm{c}_d \, \max \left\{ \Big\|
\sum_{\mathrm{I} \in [n]^d} a_{\mathrm{I}}^{\dag} \otimes
e_{\pi_{\alpha_1}(\mathrm{I}), \pi_{\beta_1}(\mathrm{I})} \Big\|_p
, \Big\| \sum_{\mathrm{I} \in [n]^d} a_{\mathrm{I}}^{\dag} \otimes
e_{\pi_{\alpha_2}(\mathrm{I}), \pi_{\beta_2}(\mathrm{I})} \Big\|_p
\right\},
\end{eqnarray*}
for some $(\alpha_1, \beta_1)$ and $(\alpha_2, \beta_2)$ in
$\mathbb{P}_d(2)$ satisfying
\begin{eqnarray*}
\mathrm{T}(\alpha_1,\beta_1) & < & \mathrm{T}(\alpha,\beta),
\\ \mathrm{T}(\alpha_2,\beta_2) & < & \mathrm{T}(\alpha,\beta).
\end{eqnarray*}
\end{claim}

\begin{remark}
\textnormal{Clearly, iteration of Claim \ref{Claim} concludes the
proof of Theorem \ref{Theorem-Degreed}.}
\end{remark}

\dem If $\mathrm{G} = \Free_{nd}$, we define $\pi: \mathrm{G}
\rightarrow \mathcal{B}(\ell_2(\mathrm{G}))$ by
\begin{eqnarray*}
\pi(g_{rs}) = \left\{ \begin{array}{ll}
\lambda_{\mathrm{G}}(g_{rs}) & \mbox{if} \ r=\mathrm{a},\mathrm{b}
\\ 1 & \mbox{otherwise} \end{array} , \right.
\end{eqnarray*}
where $\lambda_{\mathrm{G}}$ denotes the left regular
representation of $\mathrm{G}$ and $1 \le s \le n$. Clearly, the
mapping $\pi$ extends to a unitary representation of $\mathrm{G}$.
Therefore, by Fell's absorption principle we have $$\Big\|
\sum_{\mathrm{I} \in [n]^d} a_{\mathrm{I}}^{\dag} \otimes
e_{\pi_{\alpha}(\mathrm{I}), \pi_{\beta}(\mathrm{I})} \Big\|_p =
\Big\| \sum_{\mathrm{I} \in [n]^d} a_{\mathrm{I}}^{\dag} \otimes
\lambda_{\mathrm{G}} (g_{\mathrm{b}i_{\mathrm{b}}}
g_{\mathrm{a}i_{\mathrm{a}}}) \otimes e_{\pi_{\alpha}(\mathrm{I}),
\pi_{\beta}(\mathrm{I})} \Big\|_p.$$ Recall that $\alpha, \beta
\neq \emptyset$ since otherwise $(\alpha, \beta)$ would be
non-transposed. Hence we can assume that $1 \le \mathrm{b} <
\mathrm{a} \le d$. Then we define $$\begin{array}{rclcrcl}
\alpha_1 & = & \alpha \setminus \, \{\mathrm{a}\} & \qquad &
\alpha_2 & = & \alpha \cup \{\mathrm{b}\} \\ \beta_1 & = & \beta
\cup \{\mathrm{a}\} & \qquad & \beta_2 & = & \beta \setminus \,
\{\mathrm{b}\}. \end{array}$$ In particular, we can write
\begin{equation} \label{Equation-Separation}
\sum_{\mathrm{I} \in [n]^d} a_{\mathrm{I}}^{\dag} \otimes
\lambda_{\mathrm{G}} (g_{\mathrm{b}i_{\mathrm{b}}}
g_{\mathrm{a}i_{\mathrm{a}}}) \otimes e_{\pi_{\alpha}(\mathrm{I})
, \pi_{\beta}(\mathrm{I})} =
\sum_{i_{\mathrm{a}},i_{\mathrm{b}}=1}^n
x_{i_{\mathrm{b}}i_{\mathrm{a}}} \otimes \lambda_{\mathrm{G}}
(g_{\mathrm{b}i_{\mathrm{b}}} g_{\mathrm{a}i_{\mathrm{a}}})
\otimes e_{i_{\mathrm{a}},i_{\mathrm{b}}},
\end{equation}
where $x_{i_{\mathrm{b}}i_{\mathrm{a}}}$ has the following form
$$x_{i_{\mathrm{b}}i_{\mathrm{a}}} = \sum_{i_1, \ldots,
i_{\mathrm{b}-1}=1}^n \, \sum_{i_{\mathrm{b}+1}, \ldots,
i_{\mathrm{a}-1}=1}^n \, \sum_{i_{\mathrm{a}+1}, \ldots, i_d=1}^n
a_{i_1 \cdots i_d}^{\dag} \otimes e_{\pi_{\alpha_1}(i_1 \cdots
i_d), \pi_{\beta_2}(i_1 \cdots i_d)},$$ with the obvious
modifications on the sum indices if $\mathrm{a}=d$ or
$\mathrm{b}=1$ or $\mathrm{a} = \mathrm{b}+1$. On the other hand,
since $x_{i_{\mathrm{b}}i_{\mathrm{a}}}$ lives in some
non-commutative $L_p$ space $L_p(\varphi)$, we can apply Lemma
\ref{Lemma-Transposition} to the right hand side of
(\ref{Equation-Separation}) to obtain
\begin{eqnarray*}
\lefteqn{\Big\| \sum_{i_{\mathrm{a}},i_{\mathrm{b}}=1}^n
x_{i_{\mathrm{b}}i_{\mathrm{a}}} \otimes \lambda_{\mathrm{G}}
(g_{\mathrm{b}i_{\mathrm{b}}} g_{\mathrm{a}i_{\mathrm{a}}})
\otimes e_{i_{\mathrm{a}},i_{\mathrm{b}}} \Big\|_p} \\ & \le &
\mathrm{c}_d \max \left\{ \Big\|
\sum_{i_{\mathrm{a}},i_{\mathrm{b}}=1}^n
x_{i_{\mathrm{b}}i_{\mathrm{a}}} \otimes
e_{1,i_{\mathrm{b}}i_{\mathrm{a}}}
\Big\|_{L_p(\varphi;R_p^{n^2})}, \Big\|
\sum_{i_{\mathrm{a}},i_{\mathrm{b}}=1}^n
x_{i_{\mathrm{b}}i_{\mathrm{a}}} \otimes
e_{i_{\mathrm{b}}i_{\mathrm{a}},1} \Big\|_{L_p(\varphi;C_p^{n^2})}
\right\}.
\end{eqnarray*}
Finally, we observe that
\begin{eqnarray*}
\sum_{i_{\mathrm{a}},i_{\mathrm{b}}=1}^n
x_{i_{\mathrm{b}}i_{\mathrm{a}}} \otimes
e_{1,i_{\mathrm{b}}i_{\mathrm{a}}} & = & \sum_{\mathrm{I} \in
[n]^d} a_{\mathrm{I}}^{\dag} \otimes
e_{\pi_{\alpha_1}(\mathrm{I}), \pi_{\beta_1}(\mathrm{I})}, \\
\sum_{i_{\mathrm{a}},i_{\mathrm{b}}=1}^n
x_{i_{\mathrm{b}}i_{\mathrm{a}}} \otimes
e_{i_{\mathrm{b}}i_{\mathrm{a}},1} & = & \sum_{\mathrm{I} \in
[n]^d} a_{\mathrm{I}}^{\dag} \otimes
e_{\pi_{\alpha_2}(\mathrm{I}), \pi_{\beta_2}(\mathrm{I})}.
\end{eqnarray*}
This completes the proof since it is clear that
$\mathrm{T}(\alpha_1,\beta_1), \mathrm{T}(\alpha_2,\beta_2) <
\mathrm{T}(\alpha,\beta)$. \fin

\begin{remark} \label{Remark-Counterexample}
\textnormal{A consequence of our proof is the following
equivalence of norms $$\max_{0 \le k \le d} \left\{ \Big\|
\sum_{\mathrm{I} \in [n]^k} \sum_{\mathrm{J} \in[n]^{d-k}}
a_{\mathrm{IJ}}^{\dag} \otimes e_{\mathrm{I}, \mathrm{J}}^{}
\Big\|_p \right\} \simeq \max_{(\alpha,\beta) \in \mathbb{P}_d(2)}
\left\{ \Big\| \sum_{\mathrm{I} \in [n]^d} a_{\mathrm{I}}^{\dag}
\otimes e_{\pi_{\alpha}(\mathrm{I}), \pi_{\beta}(\mathrm{I})}
\Big\|_p \right\}.$$ This means that, for operators of the form
$a_{i_1 \cdots i_d} \otimes \lambda(g_{1i_1} \cdots g_{di_d})$,
the transposed terms are controlled by the non-transposed ones.
The presence of $\lambda(g_{1i_1} \cdots g_{di_d})$ is essential
to apply Fell's absorption principle. Moreover, this equivalence
is no longer true for arbitrary families of operators. A simple
counterexample is given by the 2-indexed family $a_{ij} = e_{ji}
\in S_p^n$. Namely, it is easy to check that $$\begin{array}{lcl}
\displaystyle \Big\| \sum_{i,j=1}^n e_{ji} \otimes e_{1,ij}
\Big\|_{S_p^n(R_p^{n^2})} = n^{1/2 + 1/p}, & \qquad &
\displaystyle \Big\| \sum_{i,j=1}^n e_{ji} \otimes e_{ij}
\Big\|_{S_p^n(S_p^n)} = n^{2/p}, \\ \displaystyle \Big\|
\sum_{i,j=1}^n e_{ji} \otimes e_{ij,1} \Big\|_{S_p^n(C_p^{n^2})} =
n^{1/2 + 1/p}, & \qquad & \displaystyle \Big\| \sum_{i,j=1}^n
e_{ji} \otimes e_{ji} \Big\|_{S_p^n(S_p^n)} = n.
\end{array}$$
In other words, the non-transposed term is not controlled by the
transposed ones. In particular, we conclude that the estimation
given in Section \ref{Section-Iteration} for the iteration of
Khintchine inequality is not equivalent to that provided by
Theorem \ref{Theorem-Degreed}.}
\end{remark}

\subsection{The main result for words of length $d$}
\label{Subsection-Words}

As we have pointed out several times in this paper, Buchholz's
result also holds for the whole set of reduced words of length
$d$. Therefore, it is natural to seek for the analog of Theorem
\ref{Theorem-Degreed} in this case. We shall need the following
modified version of Lemma \ref{Lemma-Transposition}.

\begin{lemma} \label{Lemma-Transposition-Inverses}
The following inequality holds for any exponent $2 \le p \le
\infty$
\begin{eqnarray*}
\lefteqn{\Big\| \sum_{i,j=1}^n a_{ij} \otimes \lambda(g_ig_j^{-1})
\otimes e_{ji} \Big\|_{L_p(\varphi \otimes \tau; S_p^n)}} \\ & \le
& \mathrm{c} \, \max \left\{ \Big\| \sum_{i,j=1}^n a_{ij} \otimes
e_{1,ij} \Big\|_{L_p(\varphi;R_p^{n^2})} , \Big\| \sum_{i,j=1}^n
a_{ij} \otimes e_{ij,1} \Big\|_{L_p(\varphi;C_p^{n^2})} \right\}.
\end{eqnarray*}
\end{lemma}

\dem We can split the sum on the left hand side as follows
$$\sum_{i,j=1}^n a_{ij} \otimes \lambda(g_ig_j^{-1}) \otimes
e_{ji} = \sum_{k=1}^n a_{kk} \otimes 1 \otimes e_{kk} + \sum_{1
\le i \neq j \le n} a_{ij} \otimes \lambda(g_ig_j^{-1}) \otimes
e_{ji}.$$ Since $g_ig_j^{-1}$ is a reduced word of length $2$
whenever $i \neq j$, the arguments employed in the proof of Lemma
\ref{Lemma-Transposition} apply to estimate the norm of the second
sum on the right in $L_p(\varphi \otimes \tau; S_p^n)$. For the
first sum, the estimation is obvious since it follows by complex
interpolation when we replace $\max$ by $\min$ above. \fin

Before stating the announced result, we redefine the analogs of
the operator spaces $\Max_p(n,d)$ and $\Word_p(n,d)$ in this new
framework. To that aim, let us define the elements $h_1, h_2,
\ldots, h_{2n}$ of $\Free_n$ as follows $$h_k = \left\{
\begin{array}{ll} g_k & \mbox{if} \ 1 \le k \le n \\ g_{k-n}^{-1}
& \mbox{otherwise} \end{array}, \right.$$ where $g_1, g_2, \ldots,
g_n$ are the generators of $\Free_n$. In this paragraph,
$\mathsf{W}_p(n,d)$ will denote the subspace of $L_p(\tau)$
generated by the image under $\lambda$ of the set of reduced words
of length $d$. In other words, an element of $\mathsf{W}_p(n,d)$
has the form $$\sum_{|t|=d} \alpha_t \lambda(t) = \sum_{i_1,
\ldots, i_d=1}^{2n} \alpha_{i_1i_2 \cdots i_d} \lambda(h_{i_1}
h_{i_2} \cdots h_{i_d}) \in L_p(\tau),$$ where the family of
scalars $$\mathsf{A} = \Big\{ \alpha_{i_1i_2 \cdots i_d} \, \Big|
\ 1 \le i_k \le 2n \Big\},$$ satisfies the following
\emph{cancellation property}
\begin{equation} \label{Equation-Cancellation}
\alpha_{i_1i_2 \cdots i_d} = 0 \quad \mbox{if} \quad i_s \equiv n
+ i_{s+1} \ (\mbox{mod} \, 2n),
\end{equation}
for some $1 \le s < d$. Note that the cancellation property
(\ref{Equation-Cancellation}) is taken so that we only consider
reduced words of length $d$. This notation will allow us to handle
the space $\mathsf{W}_p(n,d)$ just like $\Word_p(n,d)$ in the
previous paragraph. On the other hand, as pointed out in the
Introduction, we can regard the family $\mathsf{A}$ as a $(2n)^k
\times (2n)^{d-k}$ matrix as follows $$\mathsf{A}_k = \Big(
\alpha_{(i_1 \cdots i_k),(i_{k+1} \cdots i_d)} \Big) \in
C_p^{(2n)^k} \otimes_h R_p^{(2n)^{d-k}}.$$ Then, we define
$\mathsf{K}_p(n,d)$ as the subspace of $$\begin{array}{lcl}
\displaystyle \mathcal{J}_p(n,d) = \sum_{k=0}^d C_p^{(2n)^k}
\otimes_h R_p^{(2n)^{d-k}} & \qquad & \mbox{if} \ 1 \le p \le 2,
\\ \displaystyle \mathcal{J}_p(n,d) = \bigcap_{k=0}^d C_p^{(2n)^k}
\otimes_h R_p^{(2n)^{d-k}} & \qquad & \mbox{if} \ 2 \le p \le
\infty,
\end{array}$$ where certain entries are zero according to
(\ref{Equation-Cancellation}). We shall prove the following
result.

\begin{corollary} \label{Corollary-Lengthd}
$\mathsf{K}_p(n,d)$ and $\mathsf{W}_p(n,d)$ are completely
isomorphic for $1 \le p \le \infty$. More concretely, there exists
an absolute constant $\mathrm{c}_d$ depending only on $d$ such
that the following inequalities hold for any family $\mathcal{A}$
of operators in $L_p(\varphi)$ satisfying the cancellation
property $(\ref{Equation-Cancellation})$ $$\frac{1}{\mathrm{c}_d}
\, \|\mathcal{A}\|_{L_p(\varphi; \mathsf{K}_p(n,d))} \le \Big\|
\sum_{i_1 \ldots i_d=1}^{2n} a_{i_1 \cdot \cdot i_d} \otimes
\lambda(h_{i_1} \cdots  h_{i_d}) \Big\|_{L_p(\varphi \otimes
\tau)} \le \mathrm{c}_d \|\mathcal{A}\|_{L_p(\varphi;
\mathsf{K}_p(n,d))}.$$ Moreover, the natural projection
$\mathrm{P}: L_p(\tau) \rightarrow \mathsf{W}_p(n,d)$ is c.b. with
$\|\mathrm{P}\|_{cb} \le \mathrm{c}_d$.
\end{corollary}

The proof we are giving is quite similar to that of Theorem
\ref{Theorem-Degreed}. In particular, we shall skip those
arguments which already appeared above. The first remark is that
the case $p=\infty$ is exactly the content of Buchholz's result in
\cite{Bu1}. Therefore, arguing as we did after the statement of
Theorem \ref{Theorem-Degreed}, we have:
\begin{itemize}
\item The last assertion of Corollary
\ref{Corollary-Lengthd} holds.
\item The spaces $\mathsf{W}_p(n,d)$ are an
interpolation family for $1 \le p \le \infty$.
\item The dual of the space $\mathsf{W}_p(n,d)$ is
completely isomorphic to $\mathsf{W}_{p'}(n,d)$.
\item The case $1 \le p \le 2$ in Corollary
\ref{Corollary-Lengthd} follows from the case $2 \le p \le
\infty$.
\item The lower estimate for the case
$2 \le p \le \infty$ holds with some constant $\mathrm{c}_d$.
\end{itemize}

The last two points use that $\mathsf{K}_p(n,d)$ is completely
complemented in $\mathcal{J}_p(n,d)$ with constants independent on
$n$ and $p$. The proof of this fact is simple. Indeed, by
transposition and complex interpolation it suffices to prove it
for $p = \infty$. Now, since $\mathcal{J}_{\infty}(n,d)$ is an
intersection space, we just need to see it for each space
appearing in the intersection. Let $\mathcal{I}$ be the set of
indices in $[2n]^d$ satisfying the cancellation property
(\ref{Equation-Cancellation}) and let $\mathsf{H}_{\infty}(n,d)$
be the subspace of elements of $\mathcal{J}_{\infty}(n,d)$
supported in $\mathcal{I}$. Then, it is clear that the projection
$\mathrm{Q}$ onto $\mathsf{H}_{\infty}(n,d)$ is completely bounded
since it decomposes as a sum of $d-1$ \emph{diagonal} projections.
In particular, the projection $\mathrm{P}$ onto
$\mathsf{K}_p(n,d)$ is also completely bounded.

\vskip3pt

\sketch We only prove the upper estimate for $2 \le p \le \infty$.
As above, let us write $\mathrm{G}$ for the free group
$\Free_{nd}$ and $\psi_{nd}$ for the standard trace on its reduced
$C^*$-algebra. Now, following the notation just introduced, we
label the set of generators and its inverses by $h_{1k}, h_{2k},
\ldots, h_{dk}$ with $1 \le k \le 2n$. If $\lambda_{\mathrm{G}}$
denotes the left regular representation of $\mathrm{G}$, we
consider the family of operators $$\mathcal{A}^{\dag} = \Big\{
a_{i_1i_2 \cdots i_d} \otimes \lambda_{\mathrm{G}}(h_{1i_1}
h_{2i_2} \cdots h_{di_d}) \, \Big| \ 1 \le i_k \le 2n \Big\}.$$
The following chain of inequalities can be proved applying the
same arguments as for the proof of Theorem \ref{Theorem-Degreed}.
Namely, essentially we use Buchholz's result, complex
interpolation and the iteration of Khintchine inequality described
in Section \ref{Section-Iteration}.
\begin{eqnarray} \label{Equation-Chain}
\lefteqn{\Big\| \sum_{i_1, \ldots, i_d=1}^{2n} a_{i_1 \cdots i_d}
\otimes \lambda(h_{i_1} \cdots h_{i_d}) \Big\|_p} \\ \nonumber &
\le & \mathrm{c}_d \, \Big\| \sum_{i_1, \ldots, i_d=1}^{2n} a_{i_1
\cdots i_d} \otimes \lambda_{\mathrm{G}}(h_{1i_1} \cdots h_{di_d})
\Big\|_p \\ \nonumber & \le & \mathrm{c}_d \, \Big\| \sum_{i_1,
\ldots, i_d=1}^{2n} \lambda(h_{i_1}) \otimes \cdots \otimes
\lambda(h_{i_d}) \otimes a_{i_1 \cdots i_d} \otimes
\lambda_{\mathrm{G}}(h_{1i_1} \cdots h_{di_d}) \Big\|_p \\
\nonumber & \le & \mathrm{c}_d \, \max_{(\alpha,\beta) \in
\mathbb{P}_d(2)} \left\{ \Big\| \sum_{\mathrm{I} \in [2n]^d}
a_{\mathrm{I}}^{} \otimes \lambda_{\mathrm{G}}(h_{1i_1} \cdots
h_{di_d}) \otimes e_{\pi_{\alpha}(\mathrm{I}),
\pi_{\beta}(\mathrm{I})} \Big\|_p \right\}.
\end{eqnarray}
Then, the proof reduces again to the proof of
\begin{equation} \label{Equation-Plan-Lengthd}
\mathrm{RHS}_{(\ref{Equation-Chain})} \le \mathrm{c}_d \ \max_{0
\le k \le d} \left\{ \Big\| \sum_{\mathrm{I} \in [2n]^k}
\sum_{\mathrm{J} \in[2n]^{d-k}} a_{\mathrm{IJ}}^{} \otimes
\lambda_{\mathrm{G}}(h_{1i_1} \cdots h_{di_d}) \otimes
e_{\mathrm{I}, \mathrm{J}}^{} \Big\|_p \right\}.
\end{equation}
By Fell's absorption principle, we have $$\Big\| \sum_{\mathrm{I}
\in [2n]^d} a_{\mathrm{I}}^{\dag} \otimes
e_{\pi_{\alpha}(\mathrm{I}), \pi_{\beta}(\mathrm{I})} \Big\|_p =
\Big\| \sum_{\mathrm{I} \in [2n]^d} a_{\mathrm{I}}^{\dag} \otimes
\lambda_{\mathrm{G}} (h_{\mathrm{b}i_{\mathrm{b}}}
h_{\mathrm{a}i_{\mathrm{a}}}) \otimes e_{\pi_{\alpha}(\mathrm{I}),
\pi_{\beta}(\mathrm{I})} \Big\|_p,$$ with $a_{i_1 \cdots
i_d}^{\dag} = a_{i_1 \cdots i_d} \otimes
\lambda_{\mathrm{G}}(h_{1i_1} \cdots h_{di_d})$. Moreover, we can
write
\begin{equation} \label{Equation-Separation-Lengthd}
\sum_{\mathrm{I} \in [2n]^d} a_{\mathrm{I}}^{\dag} \otimes
\lambda_{\mathrm{G}} (h_{\mathrm{b}i_{\mathrm{b}}}
h_{\mathrm{a}i_{\mathrm{a}}}) \otimes e_{\pi_{\alpha}(\mathrm{I})
, \pi_{\beta}(\mathrm{I})} =
\sum_{i_{\mathrm{a}},i_{\mathrm{b}}=1}^{2n}
x_{i_{\mathrm{b}}i_{\mathrm{a}}} \otimes \lambda_{\mathrm{G}}
(h_{\mathrm{b}i_{\mathrm{b}}} h_{\mathrm{a}i_{\mathrm{a}}})
\otimes e_{i_{\mathrm{a}},i_{\mathrm{b}}},
\end{equation}
where $x_{i_{\mathrm{b}}i_{\mathrm{a}}}$ has the form
$$x_{i_{\mathrm{b}}i_{\mathrm{a}}} = \sum_{i_1, \ldots,
i_{\mathrm{b}-1}=1}^{2n} \, \sum_{i_{\mathrm{b}+1}, \ldots,
i_{\mathrm{a}-1}=1}^{2n} \, \sum_{i_{\mathrm{a}+1}, \ldots,
i_d=1}^{2n} a_{i_1 \cdots i_d}^{\dag} \otimes
e_{\pi_{\alpha_1}(i_1 \cdots i_d), \pi_{\beta_2}(i_1 \cdots
i_d)},$$ with the obvious modifications if $\mathrm{a}=d$ or
$\mathrm{b}=1$ or $\mathrm{a} = \mathrm{b}+1$. Recall that the
operators $x_{i_{\mathrm{b}},i_{\mathrm{a}}}$ do not necessarily
satisfy the cancellation property (\ref{Equation-Cancellation}).
However, we can decompose the sum in
(\ref{Equation-Separation-Lengthd}) as follows
\begin{eqnarray*}
\sum_{i_{\mathrm{a}},i_{\mathrm{b}}=1}^{2n}
x_{i_{\mathrm{b}}i_{\mathrm{a}}} \otimes \lambda_{\mathrm{G}}
(h_{\mathrm{b}i_{\mathrm{b}}} h_{\mathrm{a}i_{\mathrm{a}}})
\otimes e_{i_{\mathrm{a}},i_{\mathrm{b}}} & = & \
\sum_{i_{\mathrm{a}},i_{\mathrm{b}}=1}^{n}
x_{i_{\mathrm{b}}i_{\mathrm{a}}} \otimes \lambda_{\mathrm{G}}
(h_{\mathrm{b}i_{\mathrm{b}}} h_{\mathrm{a}i_{\mathrm{a}}})
\otimes e_{i_{\mathrm{a}},i_{\mathrm{b}}} \\ & + &
\sum_{i_{\mathrm{a}}=1}^{n} \, \sum_{i_{\mathrm{b}}=n+1}^{2n}
x_{i_{\mathrm{b}}i_{\mathrm{a}}} \otimes \lambda_{\mathrm{G}}
(h_{\mathrm{b}i_{\mathrm{b}}} h_{\mathrm{a}i_{\mathrm{a}}})
\otimes e_{i_{\mathrm{a}},i_{\mathrm{b}}} \\ & + &
\sum_{i_{\mathrm{b}}=1}^{n} \, \sum_{i_{\mathrm{a}}=n+1}^{2n}
x_{i_{\mathrm{b}}i_{\mathrm{a}}} \otimes \lambda_{\mathrm{G}}
(h_{\mathrm{b}i_{\mathrm{b}}} h_{\mathrm{a}i_{\mathrm{a}}})
\otimes e_{i_{\mathrm{a}},i_{\mathrm{b}}} \\ & + &
\sum_{i_{\mathrm{a}},i_{\mathrm{b}}=n+1}^{2n}
x_{i_{\mathrm{b}}i_{\mathrm{a}}} \otimes \lambda_{\mathrm{G}}
(h_{\mathrm{b}i_{\mathrm{b}}} h_{\mathrm{a}i_{\mathrm{a}}})
\otimes e_{i_{\mathrm{a}},i_{\mathrm{b}}}.
\end{eqnarray*}
Then it is clear than Lemma \ref{Lemma-Transposition} applies to
the first and the fourth sums while Lemma
\ref{Lemma-Transposition-Inverses} applies to the second and third
sums. In summary, we have
\begin{eqnarray*}
\lefteqn{\Big\| \sum_{i_{\mathrm{a}},i_{\mathrm{b}}=1}^{2n}
x_{i_{\mathrm{b}}i_{\mathrm{a}}} \otimes \lambda_{\mathrm{G}}
(h_{\mathrm{b}i_{\mathrm{b}}} h_{\mathrm{a}i_{\mathrm{a}}})
\otimes e_{i_{\mathrm{a}},i_{\mathrm{b}}} \Big\|_p} \\ & \le &
\mathrm{c}_d \max \left\{ \Big\|
\sum_{i_{\mathrm{a}},i_{\mathrm{b}}=1}^{2n}
x_{i_{\mathrm{b}}i_{\mathrm{a}}} \otimes
e_{1,i_{\mathrm{b}}i_{\mathrm{a}}}
\Big\|_{L_p(\varphi;R_p^{n^2})}, \Big\|
\sum_{i_{\mathrm{a}},i_{\mathrm{b}}=1}^{2n}
x_{i_{\mathrm{b}}i_{\mathrm{a}}} \otimes
e_{i_{\mathrm{b}}i_{\mathrm{a}},1} \Big\|_{L_p(\varphi;C_p^{n^2})}
\right\}.
\end{eqnarray*}
Finally, we conclude as in Claim \ref{Claim}. This completes the
proof of (\ref{Equation-Plan-Lengthd}). \fin

\begin{remark}
\textnormal{In Voiculescu's free probability theory, stochastic
independence of random variables is replaced by freeness of
non-commutative random variables. In this setting, the Wigner's
probability distribution $$d\mu_{\mathrm{W}}(t) = 1_{[-2,2]}
\frac{\sqrt{4 - t^2}}{2\pi} \, dt$$ plays a crucial role. Namely,
given a free family $x_1, x_2, \ldots, x_n$ of self-adjoint random
variables in a non-commutative probability space $(\mathcal{M},
\tau)$, we say that $x_1, x_2, \ldots, x_n$ is a \emph{free
semi-circular system} if each $x_k$ is equipped with Wigner's
distribution. This family is the free analog of a system of $n$
independent standard real-valued gaussian random variables.
Explicit constructions of free semi-circular systems are available
by means of the creation and annihilation operators on the full
Fock space, see \cite{P2,VDN} for more on this. The free analog of
$n$ independent complex-valued gaussians is now given by taking
$$z_k = \frac{1}{\sqrt{2}} (x_k' + i x_k''),$$ with $x_1', x_1'',
x_2', x_2'', \ldots, x_n', x_n''$ being a free semi-circular
system. This new system is called a \emph{free circular system}.
At this point, it is natural to guess that the analog of Theorem
\ref{Theorem-Degreed} should hold when we replace free generators
by free circular random variables. Indeed, as it was pointed out
in the Introduction, the family of operators $\lambda(g_1),
\lambda(g_2), \ldots, \lambda(g_n)$ is the free analog of the
sequence of Rademacher functions $\rad_1, \rad_2, \ldots, \rad_n$.
Therefore, a free version of the central limit theorem is exactly
what is needed here. A precise statement of this result can be
found in \cite{VDN} and supports the previous identification
between real-valued gaussians and semi-circular random variables.
Although we are not giving the details, it can be checked that the
central limit theorem for free random variables provides the
analog of Theorem \ref{Theorem-Degreed} for free circular
variables. In other words, if we replace the operators
$\lambda(g_{i_1}g_{i_2} \cdots g_{i_d})$ by the products $z_{i_1}
z_{i_2} \cdots z_{i_d}$ in Theorem \ref{Theorem-Degreed}, then the
same conclusions hold. In passing, we also refer the interested
reader to Nou's paper \cite{N}, which contains the analog of
Buchholz's result for $q$-gaussian randon variables.}
\end{remark}

\begin{remark}
\textnormal{The paper \cite{P1} deals with the notion of
$p$-orthogonal sums in non commutative $L_p$ spaces. Applying some
combinatorial techniques, it is shown that the Khintchine type
inequality that applies for $\Word_p(n)$ majorizes the behaviour
of a much larger class of operators, the so-called $p$-orthogonal
sums, for any even integer $p$. On the other hand, the bounds
given in Section \ref{Section-Iteration} for $\Word(n)^{\otimes
d}$ constitute an upper bound of a more general family of
operators. Namely, let $(\mathcal{M},\tau)$ be a von Neumann
algebra endowed with a standard trace satisfying $\tau(1)=1$ and
let $L_p(\tau)$ be the associated non-commutative $L_p$ space. Let
$\Gamma$ stand for the product set $[n] \times \cdots \times [n]$
with $d$ factors. Then, given an even integer $p$ and a family $f
= (f_\gamma)_{\gamma \in \Gamma}$ of operators in $L_p(\tau)$
indexed by $\Gamma$, we shall say that $f$ is $p$-orthogonal with
$d$ indices if $$\tau \big(
f_{h(1)}^*f_{h(2)}^{}f_{h(3)}^*f_{h(4)}^{} \cdots
f_{h(p-1)}^*f_{h(p)}^{} \big) = 0$$ whenever the function $h:
\{1,2, \ldots, p\} \rightarrow \Gamma$ has an injective
projection. In other words, whenever the coordinate function
$\pi_k \circ h: \{1,2, \ldots, p\} \rightarrow [n]$ is an
injective function for some index $1 \le k \le d$. The paper
\cite{Pa} extends the results in \cite{P1} to this more general
setting by studying the norm in $L_p(\tau)$ of the sum
$$\sum_{\gamma \in \Gamma} f_{\gamma}.$$ More concretely, the norm
of this sum in $L_p(\tau)$ is bounded above by the expressions
given in Section \ref{Section-Iteration}, see \cite{Pa} for a
precise statement. Moreover, we should point out that, in contrast
with the image under $\lambda$ of the words of length $d$, the
family $\lambda(g_{i_1}g_{i_2} \cdots g_{i_d})$ is also a
$p$-orthogonal family with $d$ indices.}
\end{remark}

\noindent \textbf{Acknowledgements.} The first author was
partially supported by the Project MTM2004-00678, Spain. The
second author was partially supported by the NSF and by the Texas
Advanced Research Program 010366-163. This work was carried out
while the first-named author was a Visiting Asistant Professor at
Texas A\&M University. The first-named author would like to thank
the Math Department for its support and hospitality.

\

\noindent \textbf{Note added in proof.}  Recently, Ricard and Xu
\cite{RX} have extended Buchholz's result \cite{Bu1} to arbitrary
free product $C^*$-algebras. As they point out in their paper, the
same construction holds for amalgamated free products of von
Neumann algebras. Moreover, after Ricard/Xu's work, Junge and the
first-named author have generalized in \cite{JP} the main result
in \cite{RX} to arbitrary indices $2 \le p \le \infty$ as a
consequence of the free analogue of Rosenthal's inequality
\cite{Ro}, also proved in \cite{JP}.

\bibliographystyle{amsplain}

\begin{thebibliography}{10}
\bibitem {Bu1} A. Buchholz, \emph{Norm of convolution by
operator-valued functions on free groups}, Proc. Amer. Math. Soc.
\textbf{127} ($1999$), $1671-1682$.
\bibitem {Bu2} A. Buchholz, \emph{Operator Khintchine inequality
in non-commutative probability}, Math. Ann. \textbf{319} ($2001$),
$1-16$.
\bibitem {HP} U. Haagerup and G. Pisier, \emph{Bounded linear
operators between $C^*$-algebras}, Duke Math. J. \textbf{71}
($1993$), $889-925$.
\bibitem {J1} Personal communication of M. Junge.
\bibitem {J2} M. Junge, \emph{Embedding of the operator space OH
and the logarithmic `little Grothendieck inequality'}. To appear
in Invent. Math.
\bibitem {JP} M. Junge and J. Parcet, \emph{Rosenthal type
inequalities for amalgamated free products and applications}.
Preprint.
\bibitem {L} F. Lust-Piquard, \emph{In\'{e}galit\'{e}s de Khintchine dans
$C_p$ $(1 < p < \infty)$}, C.R. Acad. Sci. Paris \textbf{303}
($1986$), 289-292.
\bibitem {LP} F. Lust-Piquard and G. Pisier, \emph{Non-commutative
Khintchine and Paley inequalities}, \emph{Ark. Mat.} \textbf{29}
($1991$), 241-260.
\bibitem{N} A. Nou, \emph{Non injectivity of the $q$-deformed von
Neumann algebra}. Math. Ann. \textbf{330} (2004), 17-38.
\bibitem {Pa} J. Parcet, \emph{Multi-indexed $p$-orthogonal sums
in non-commutative Lebesgue spaces}. Indiana Univ. Math. J.
\textbf{53} (2004), 1171-1188.
\bibitem {P1} G. Pisier, \emph{An inequality for $p$-orthogonal
sums in non-commutative $L_p$}, Illinois J. Math. \textbf{44}
($2000$), $901-923$.
\bibitem {P2} G. Pisier, \emph{Introduction to
Operator Space Theory}, Cambridge Univ. Press, 2003.
%\bibitem {PX} G. Pisier and Q.Xu, \emph{Non-commutative martingale
%inequalities}, Comm. Math. Physics \textbf{189} (1997), $667-698$.
\bibitem {RX} E. Ricard and Q. Xu, \emph{Khitnchine type
inequalities for reduced free products and applications}.
Preprint.
\bibitem {Ro} H.P. Rosenthal, \emph{On the subspaces of
$L^p$ $(p>2)$ spanned by sequences of independent random
variables}. Israel J. Math. \textbf{8} (1970), 273-303.
\bibitem {VDN} D. Voiculescu, K. Dykema and A. Nica,
\emph{Free random variables}, CRM Monograph Series \textbf{1},
Amer. Math. Soc., 1992.
\end{thebibliography}

\

\noindent Departamento de Matem{\'a}ticas \\ Universidad Aut\'{o}noma
de Madrid \\ 28049 Madrid, Spain \\ E-mail:
\texttt{javier.parcet@uam.es}

\

\noindent \'Equipe d'Analyse \\ Universit{\'e} Paris VI \\ Case 186,
F-75252 Paris Cedex 05, France \\ and \\ Math Department \\ Texas
A\&M University \\ College Station, TX 77843, USA
\\ E-mail: \texttt{gip@ccr.jussieu.fr}

\

\noindent \textsc{Key words and phrases}: Non-commutative
Khintchine inequality, Fell's absorption principle, free group.

\vskip5pt

\noindent \textsc{2000 Mathematics Subject Classification}: 46L07,
46L52, 46L54
\end{document}